\documentclass[11pt,reqno]{amsart}
\usepackage{xcolor}
\usepackage[T1]{fontenc}
\usepackage{lmodern}
\usepackage{amsmath,amssymb,amsthm,mathtools}
\usepackage[margin=1in]{geometry}
\usepackage[hidelinks]{hyperref}
\usepackage{fourier}

\usepackage[backend=bibtex,
    sorting=anyt,
    isbn=false,
    url=false,
    doi=false,
    maxbibnames=100
    ]{biblatex}
\theoremstyle{plain}
\newtheorem*{itheorem}{Theorem}
\newtheorem{conjecture}{Conjecture}
\newtheorem{theorem}{Theorem}[section]
\newtheorem{lemma}[theorem]{Lemma}
\newtheorem{proposition}[theorem]{Proposition}
\newtheorem{corollary}[theorem]{Corollary}
\theoremstyle{definition}
\newtheorem{definition}[theorem]{Definition}

\theoremstyle{remark}

\newcommand{\bk}{\Bbbk}

\DeclareMathOperator{\CDi}{\mathsf{CD}}
\DeclareMathOperator{\Nov}{\mathsf{Nov}}
\DeclareMathOperator{\Lie}{\mathsf{Lie}}
\DeclareMathOperator{\WLie}{\mathsf{WLie}}
\DeclareMathOperator{\MC}{MC}
\DeclareMathOperator{\Tw}{Tw}
\DeclareMathOperator{\Hom}{Hom}
\DeclareMathOperator{\End}{End}
\DeclareMathOperator{\Ind}{Ind}
\DeclareMathOperator{\Res}{Res}

\newcommand{\sgn}{\mathrm{sgn}}

\newcommand{\cD}{\mathcal D}
\newcommand{\cH}{\mathcal H}

\DeclareMathAlphabet{\pazocal}{OMS}{zplm}{m}{n}

\def\calH{\pazocal{H}}

\def\calM{\pazocal{M}}

\def\calP{\pazocal{P}}

\def\calV{\pazocal{V}}

\title{Invariant operations on vector fields on the line are Lie polynomials}
\author{Vladimir Dotsenko}
\address{Institute of Mathematics, Academy of Mathematics and Systems Science, Chinese Academy of Sciences, Beijing, 100190, China}
\address{Institut de Recherche Math\'ematique Avanc\'ee, UMR 7501, Universit\'e de Strasbourg et CNRS, 7 rue Ren\'e-Descartes, 67000 Strasbourg, France}

\email{vdotsenko@unistra.fr}

\author{Farukh Mashurov}
\address{SICM, Southern University of Science and Technology, Shenzhen, 518055, China}
\address{SDU University, Kaskelen, Kazakhstan}

\email{f.mashurov@gmail.com}

\author{Bauyrzhan Sartayev}
\address{Narxoz University, Almaty, Kazakhstan}
\address{Institute of Mathematics and Mathematical Modeling, Almaty, Kazakhstan}

\email{baurjai@gmail.com}


\begin{document}

\begin{abstract}
We prove that all $W_1$-equivariant multilinear operations on the Witt Lie algebra $W_1$ are Lie polynomials. We study this question in terms of Wronskian Lie brackets and Novikov products, which are abstractions of the Lie bracket $[a,b]:=ab'-ba'$ and the nonassociative product $a\circ b=ab'$ one can define on a commutative differential algebra. In terms of these notions, our main result asserts that the operad of Wronskian Lie algebras is isomorphic to the operad of equivariant multilinear operations on~$W_1$; we also prove a parallel result for the operad of Novikov algebras, characterizing it as the operad of affine equivariant multilinear operations on $W_1$. Two immediate consequences of our work are a characterization of Lie polynomials in free Novikov algebras and a computation of the degree zero homology of the operadic twisting of the Novikov operad. Finally, while it is known that over a field of zero characteristic, a Novikov algebra embeds into its differential envelope, we show that the same is not necessarily true for Wronskian Lie algebras.  
\end{abstract}

\maketitle

\section{Introduction}

It is well known that if $A$ is a commutative differential algebra, then the operation $[a,b]:=ab'-ba'$ defines a Lie algebra structure on $A$.  This operation has appeared in the literature prominently, for instance, in work of Herz \cite{MR69810} concerning the algebraic study of PDEs, in work of Kirillov \cite{MR438390} and Lichnerowicz \cite{MR524629} on local Lie algebras and Jacobi manifolds, in work of Gelfand and Dorfman \cite{GD} on algebraic aspects of the Hamiltonian formalism, etc. Interestingly enough, it also appears in the work of Lawvere \cite{MR231882} as an example of a functor whose adjoint is worth investigating; he called this construction the Wronskian functor and asked whether a Lie algebra always embeds in its Wronskian envelope. Following the existing terminology for Lie and Jordan subalgebras of associative algebras, it is natural to call a Lie algebra that embeds in its Wronskian envelope \emph{Wronskian special}.

Notably, if $A=\bk[x]$ and $a'(x)$ is the derivative with respect to $x$, the above Lie bracket corresponds to the Lie bracket of vector fields, once we assign to an element $a(x)$ the vector field $a(x)\partial$; this Lie algebra is often denoted $W_1$ and is referred to as the Witt algebra of vector fields on the affine line. In the 1970s, several researchers including Bergman \cite{bergman1978}, Razmyslov \cite{MR794957}, and Sumenkov \cite{Sumenkov1976} discovered independently that the Lie algebra $W_1$ satisfies the identity
 \[
\sum_{\sigma \in S_{4}} \sgn(\sigma)
  \bigl[a_{\sigma(1)}, [a_{\sigma(2)}, [a_{\sigma(3)},
  [a_{\sigma(4)}, a_{5}]]]\bigr] = 0,
 \]
known as the standard Lie identity of degree five. The state-of-the-art concerning identities of the Lie algebra $W_1$ is truly mysterious:  Molev~\cite[Th.~1]{Mo} determined the dimensions and the module structures over the symmetric groups of the spaces of multilinear Lie polynomials on $W_1$, and yet we do not know the complete set of identities of this Lie algebra. It was conjectured in the 1980s \cite{MR1099432,MR794957} that all such identities follow from the standard identity of degree five, and this still remains an open question \cite[Sec.~14.3.2]{MR4368864}.  

In fact, the standard Lie identity of degree five is easily seen to hold for any Wronskian Lie bracket, and instantly shows that most Lie algebras are not Wronskian special. In recent years, Poinsot \cite{Poinsot2013,MR3424596,MR3780535} produced several papers developing the basics of Wronskian envelopes of Lie algebras; he pursued the categorical viewpoint of Lawvere and made use of the existing body of work on identities in Lie algebras. Notably, in \cite[Th.~16]{Poinsot2013} he used results of Razmyslov \cite{MR794957} to prove that subalgebras of direct sums of Lie algebras where each summand is either abelian or simple and satisfies the standard identity of degree five are Wronskian special. Later, Pogudin and Razmyslov \cite{MR3550862} established that semiprime Lie algebras satisfying the standard identity of degree five are Wronskian special. 

Another natural algebraic structure one can construct out of a commutative differential algebra is the ``half of the Wronskian bracket'', that is, the product $a\circ b=ab'$. For this product, it was proved by Dzhumadil'daev and Löfwall \cite{DL} that all identities satisfied by it follow from the two identities 
\begin{gather*}
(a\circ b)\circ c-a\circ (b\circ c)=(b\circ a)\circ c-b\circ (a\circ c),\\
(a\circ b)\circ c=(a\circ c)\circ b.
\end{gather*}
Algebras satisfying these two identities are called \emph{Novikov algebras}; they first appeared in the work of Gelfand and Dorfman \cite[\S 6]{GD}. They were given their name by Osborn \cite{Os}, presumably to highlight Novikov's idea to use these algebras in the study of Poisson brackets of hydrodynamic type~\cite{BN}.

There are two particularly important structural results concerning Novikov algebras. Bokut, Chen and Zhang \cite{MR3590862} proved that one can faithfully represent Novikov algebras in their commutative differential envelopes: every Novikov algebra $N$ over a field of zero characteristic embeds into its commutative differential envelope where the product becomes precisely $a\circ b=ab'$. Makar-Limanov and Umirbaev \cite{MR2884985} proved that the commutative differential algebra $\bk[x]$, viewed as a Novikov algebra, does not satisfy any nontrivial Novikov identities. A counterpart of these two properties for Wronskian Lie algebras stating that the class of Wronskian special Lie algebras coincides with the variety of Lie algebras generated by the Witt Lie algebra $W_1$ was recently conjectured by Tulenbaev, Umirbaev, and Zhelyabin \cite[Conj.~2]{MR4556333}. In this paper, we show that Wronskian special Lie algebras do not form a variety, thus disproving that conjecture. Concretely, we prove

\begin{itheorem}[{Th.~\ref{thm:exceptional}}]
For the free two-generated Wronskian Lie algebra $\WLie(x,y)$, and its elements
 \[
r=[[x,y],[[x,y],y]] \quad \text{ and } \quad  s=[[[x,y],[[x,y],x]],x],   
 \]
the Lie algebra $\WLie(x,y)/(r,s)$ is not Wronskian special. 
\end{itheorem}

Nevertheless, it is perfectly legitimate to consider the variety of Lie algebras that satisfy all the universal identities of the Lie bracket $ab'-ba'$ on a commutative differential algebra; this is precisely what we mean by Wronskian Lie algebras. We denote by $\WLie$ the operad encoding this class of algebras. In this paper, we give several new equivalent descriptions of the operad $\WLie$. Our main result is   

\begin{itheorem}[{Prop.~\ref{prop:suboperads}(2), Th.~\ref{th:LieNovInv}}]
The operad of Wronskian Lie algebras is isomorphic to the operad of multilinear Lie polynomials on $W_1$, which in turn is isomorphic to the operad of all $W_1$-equivariant multilinear operations on $W_1$. 
\end{itheorem}

We also prove (in Prop.~\ref{prop:suboperads}(1) and Corollary \ref{cor:b-Nov}) a parallel result for the operad of Novikov algebras, characterizing it as the operad of multilinear Novikov polynomials on $W_1$ and as the operad of affine equivariant multilinear operations on $W_1$. Together with the recent work of Chapoton and the first-named author \cite{MR5093214} concerning the operad of equivariant multilinear operations on the adjoint module of the Lie algebra $\mathfrak{sl}_2$, these results suggest that it would be interesting to study operads of equivariant multilinear operations for other Lie algebras and their modules as well. 

As an easy consequence of our work, we were able to answer the question of Dzhumadildaev and Ismailov \cite{DI14} on recognizing Lie polynomials in free Novikov algebras, proving the following criterion. 

\begin{itheorem}[{Th.~\ref{th:criterion-main}}]
Let $X$ be a finite set, and let $y\notin X$. Denote by $\delta_y$ the derivation of the free Novikov algebra $\Nov(X\cup\{y\})$ defined on the generators by $\delta_y(x_i)=[y,x_i]$ and $\delta_y(y)=0$. An element  $f\in\Nov(X)$ is a Lie polynomial if and only if
\begin{equation*}
\delta_y(f)=[y,f].    
\end{equation*}
\end{itheorem}

As a by-product of the latter result, we were able to exhibit yet another interpretation of the operad of Wronskian Lie algebras in Corollary \ref{cor:twisting-nov}: it is isomorphic to the degree zero homology of the operadic twisting \cite{MR4621635} of the operad of Novikov algebras.

\subsection*{Acknowledgements} This work was supported by the Science Committee of the Ministry of Science and Higher Education of the Republic of Kazakhstan (Grant No. BR 28713025). The first-named author is indebted to Alexander Molev for discussions of the results of the paper and very helpful historical comments. The final draft of the paper was prepared during the visit of the first-named author and the third-named author to Shenzhen International Center for Mathematics, and they thank that institution for its hospitality and excellent working conditions.

\subsection*{Disclosure of the AI use}
The present work is written by humans and for humans. In a prior preprint of the second-named and the third-named authors \cite{MS2026}, proofs of two of the statements of Corollary~\ref{cor:criterion-symbols}, obtained with assistance of ChatGPT Astra, were recorded. Once that preprint was disseminated, it became apparent that the methods used in it could be combined with unpublished work of the first-named author to prove the main theorem of the paper and to obtain a more effective criterion for Lie polynomials, and we started a collaborative project that resulted in the present paper. It contains a substantial body of new results; moreover, the proofs of results that were already present in \cite{MS2026} were completely rewritten to make them human-accessible, some statements used in those proofs were replaced by either simpler ones or more general ones, and provenance of various logical steps, not at all clear from the AI-generated text, was traced, so that credit could be given where credit is due. 

\subsection*{Conventions and notation}

All vector spaces in this paper are defined over a field $\bk$ of zero characteristic. We use the convention $\mathbb{N}=\mathbb{Z}_{\ge 0}$.  We denote by $W_1$ the Lie algebra of polynomial vector fields on the affine line; it is the vector space with the basis $X_i=x^{i+1}\partial$, $i\ge -1$, with the Lie bracket given by $[X_i,X_j]=(j-i)X_{i+j}$. For any $\alpha\in\bk$, this Lie algebra has a module called the module of densities of weight $\alpha$ and denoted $F_\alpha$; its elements are formal expressions $g(x)dx^\alpha$ with $g(x)\in\bk[x]$, with the action given by
 \[
f(x)\partial(g(x)dx^\alpha)=\left(f(x)g'(x)+\alpha f'(x)g(x)\right) dx^\alpha.
 \]
In particular, for $\alpha=-1$ this is the adjoint representation of $W_1$. 

We switch quite freely between the language of PI-algebras and the more categorical language of the operad theory. Our choice is to not overload this short article by excessive recollections, referring the reader to the monographs \cite{MR3642294,MR2954392} for basics of operad theory, and to \cite[Sec.~2.3]{MR4675074} for introduction to the language of operads intended specifically for experts in PI-theory. In particular, for a vector space $V$, we denote by $\End_V$ the endomorphism operad of $V$ consisting of all multilinear maps on $V$:
 \[
\End_V(n)=\Hom(V^{\otimes n},V),     
 \]
and if $V$ happens to be a module for a Lie algebra $\mathfrak{g}$, we denote by $\End_V^{\mathfrak{g}}$ the equivariant endomorphism operad of $V$: 
 \[
\End_V^{\mathfrak{g}}(n)=\Hom_{\mathfrak{g}}(V^{\otimes n},V).     
 \]

\section{Commutative differential algebras, Novikov algebras, and Wronskian Lie algebras}\label{sec:CD-Nov-Wr}

A commutative differential algebra \cite{Kolchin,Ritt} is a commutative associative algebra equipped with a derivation; we shall denote the product by $ab$, the derivation by $a'$, and the $k$-fold iteration of the derivation by $a^{(k)}$. We shall denote by $\CDi$ the operad of commutative differential algebras. Elements of the components $\CDi(n)$ are linear combinations of multilinear differential monomials $a_1^{(d_1)}\cdots a_n^{(d_n)}$, where the sum $d_1+\cdots+d_n$ is called the \emph{differential order of the monomial}; every universal multilinear operation one can define on a commutative associative algebra with a derivation is a linear combination of operations of this form.

It follows from results of Dzhumadildaev and Löfwall \cite{DL} that the suboperad of the operad $\CDi$ generated by the operation $a\circ b=ab'$ is spanned, in each arity $n$, by multilinear differential monomials of differential order exactly $n-1$, and that this suboperad is isomorphic to the operad of Novikov algebras. Let us record an analogous definition for the operad of Wronskian Lie algebras.

\begin{definition}
The \emph{operad of Wronskian Lie algebras} $\WLie$ is the suboperad of the operad $\CDi$ generated by the operation $[a,b]=ab'-ba'$. 
\end{definition}

This definition does not immediately furnish much information; as we discussed in the introduction, the operation $[a,b]$ is a Lie bracket, and it satisfies the standard Lie identity of degree five, but beyond that the definition does not say much. Let us record a simple result that leads to a slightly better understanding of the operad $\WLie$. 

\begin{proposition}\label{prop:suboperads}\leavevmode
\begin{enumerate}
\item The suboperad of the endomorphism operad $\End_{W_1}$ generated by the binary operation 
 \[
 f(x)\partial\circ g(x)\partial:=f(x)g'(x)\partial    
 \]
is isomorphic to the Novikov operad. 
\item The suboperad of the endomorphism operad $\End_{W_1}$ generated by the binary operation 
 \[
[f(x)\partial,g(x)\partial]:=(f(x)g'(x)-g(x)f'(x))\partial     
 \]
is isomorphic to the operad of Wronskian Lie algebras.
\end{enumerate}
\end{proposition}

\begin{proof}
The first statement is a reformulation of a result of Makar-Limanov and Umirbaev \cite{MR2884985} mentioned in the introduction who proved that the commutative differential algebra $\bk[x]$, viewed as a Novikov algebra, does not satisfy any nontrivial Novikov identities. Under the vector space isomorphism $\bk[x]\cong W_1$, the Novikov product on the polynomial algebra becomes $f(x)\partial\circ g(x)\partial:=f(x)g'(x)\partial$, completing the proof. 
The second statement is an immediate consequence of the first one: we have 
 \[
[f(x)\partial,g(x)\partial]=f(x)\partial\circ g(x)\partial-g(x)\partial\circ f(x)\partial,     
 \]
so all Lie identities of the operation $[f(x)\partial,g(x)\partial]$ are precisely all consequences of Novikov identities of the operation $f(x)\partial\circ g(x)\partial$.
\end{proof}

The result we just proved, paired with work of Molev \cite{Mo}, allows us to describe the symmetric group actions on the components of the operad $\WLie$. To state the result, let us recall the notation for the classical $S_n$-coinvariant ring 
 \[
\cH_n:=\bk[t_1,\ldots,t_n]/(\bk[t_1,\ldots,t_n]_+^{S_n}).     
 \]
Since $\cH_n$ is a quotient by a homogeneous ideal, it inherits the degree function from the polynomial ring: we denote by $\cH_n^d$ the homogeneous component of degree $d$.

\begin{corollary}\label{cor:Molev-WLie}
For all $n\ge 1$, we have an isomorphism of $S_n$-modules 
 \[
\WLie(n)\cong\cH_n^{n-1}.     
 \]
In particular, the dimension of $\WLie(n)$ is equal to the number of permutations in $S_n$ with $n-1$ inversions.
\end{corollary}

\begin{proof}
By Proposition \ref{prop:suboperads}, the operad $\WLie$ is isomorphic to the suboperad of the operad $\End_{W_1}$ generated by the binary operation $[f(x)\partial,g(x)\partial]:=(f(x)g'(x)-g(x)f'(x))\partial$.     
The latter is precisely what is studied in \cite[Th.~1]{Mo}, and the corollary follows (the statement that the dimension of $\cH_n^k$ is equal to the number of permutations in $S_n$ with $k$ inversions is well known). 
\end{proof}

The second part of Proposition \ref{prop:suboperads} is based on the fact that the Wronskian functor from commutative differential algebras to Lie algebras is a composition of two functors: the functor sending a commutative differential algebra $A$ to the Novikov algebra on the same vector space with the operation $a\circ b=ab'$, and the functor sending a Novikov algebra $A$ to a Lie algebra on the same vector space with the Lie bracket $[a,b]=a\circ b-b\circ a$. The notion of a \emph{Wronskian special} Lie algebra mentioned in the introduction concerns the left adjoint of the composite functor, but the latter functor alone also admits a meaningful left adjoint, and we may call a Lie algebra \emph{Novikov special} if it embeds into its universal enveloping Novikov algebra. 

\begin{proposition}\label{prop:NovWroSpec}
A Lie algebra $L$ is Novikov special if and only if it is Wronskian special.
\end{proposition}

\begin{proof}
This is an immediate consequence of work of Bokut, Chen, and Zhang \cite{MR3590862} who established that over a field of zero characteristic every Novikov algebra embeds into its universal enveloping commutative differential algebra. Indeed, the left adjoint functor from Lie algebras to commutative differential algebras is the composition of two left adjoints: from Lie algebras to Novikov algebras and from Novikov algebras to commutative differential algebras. Since the second one embeds every algebra into its universal envelope \cite[Th.~3]{MR3590862}, the only obstruction to embedding is given by the first one, proving the required assertion.
\end{proof}

We shall now demonstrate that a homomorphic image of a Wronskian special Lie algebra does not have to be Wronskian special, so that Wronskian special Lie algebras do not form a variety. This disproves the expectation of Tulenbaev, Umirbaev, and Zhelyabin \cite[Conj.~2]{MR4556333} that the class of Wronskian special Lie algebras coincides with the variety of Lie algebras generated by $W_1$.

\begin{theorem}\label{thm:exceptional}
For the free two-generated Wronskian Lie algebra $\WLie(x,y)$, and its elements
 \[
r=[[x,y],[[x,y],y]] \quad \text{ and } \quad  s=[[[x,y],[[x,y],x]],x],   
 \]
the Lie algebra $\WLie(x,y)/(r,s)$ is not Wronskian special.
\end{theorem}

\begin{proof}
Using Proposition \ref{prop:NovWroSpec}, we see that we do not need to distinguish between the notions of Wronskian special and Novikov special Lie algebras. Let us view $\WLie(x,y)$ as a Lie subalgebra of the free two-generated Novikov algebra $\Nov(x,y)$. 
We shall prove that the element 
\[
 w=[[x,y],[[[x,y],x],[[x,y],y]]]
\]
does not belong to the ideal of $\WLie(x,y)$ generated by $r$ and $s$, but belongs to the ideal of $\Nov(x,y)$ generated by those elements. This will establish the result: the coset of the element $w$ is a non-zero element of the algebra $\WLie(x,y)/(r,s)$ which vanishes in the universal Novikov envelope of that algebra. All computational claims below can be verified by hand; a \texttt{Python} script that performs the relevant computation is available in the online addendum \cite{addendum}.

To demonstrate that $w$ is contained in the ideal of $\Nov(x,y)$ generated by $r$ and $s$, we simply exhibit a formula for it that proves this containment:
\begin{equation}\label{id: w in J and L}
 w=([r,y]\circ x)\circ x+[x\circ r,x\circ y] +[y,s\circ y]-2(y\circ s)\circ y.
\end{equation} 

To prove that the element $w$ does not belong to the ideal of $\WLie(x,y)$ generated by $r$ and $s$, note that if we consider the bi-degree counting occurrences of $x$ and $y$, the elements $r$, $s$, $w$ are bi-homogeneous, with $\deg r=(2,3)$, $\deg s=(4,2)$, $\deg w=(4,4)$. 
The ideal of $\WLie(x,y)$ generated by $r$ and~$s$ is spanned by left-normed Lie brackets of $r,s$  with the generators $x,y$. That ideal is bi-homogeneous, and its component of multidegree $(4,4)$ is spanned by
$[[[r,y],x],x]$, $[[[r,x],y],x]$, $[[[r,x],x],y]$, and $[[s,y],y]$.
Thus, if $w$ belongs to this ideal, there exist scalars $\alpha_1,\alpha_2,\alpha_3,\alpha_4\in \bk$ such that
\begin{equation}\label{eq:linear-combination}
\begin{aligned}
 w=&\alpha_1[[[r,y],x],x]+\alpha_2[[[r,x],y],x]+\alpha_3[[[r,x],x],y]+\alpha_4[[s,y],y].
\end{aligned}
\end{equation}
Let us expand everything in the free differential algebra and compute the coefficients of the differential monomials 
\[
M_1=x^4y(y')^2y^{(5)},\qquad
M_2=x^4(y')^3y^{(4)}, \qquad
M_3=x^2x''x'''y^2(y')^2
\]
in the expansion of each element of interest. These coefficients are as follows:
\[
\renewcommand{\arraystretch}{1.2}
\begin{array}{c|rrr}
 &M_1&M_2&M_3\\ \hline
 \relax [[[r,y],x],x]&4&2&6\\
 \relax [[[r,x],y],x]&3&0&10\\
 \relax [[[r,x],x],y]&2&-2&14\\
 \relax [[s,y],y]&-1&1&-7\\
 w&0&1&-1
\end{array}
\]
Comparing the coefficients of $M_1,M_2,M_3$ on both sides of \eqref{eq:linear-combination}, we obtain
 \[
\begin{cases}
4\alpha_1+3\alpha_2+2\alpha_3-\alpha_4=0,\\
2\alpha_1-2\alpha_3+\alpha_4=1,\\
6\alpha_1+10\alpha_2+14\alpha_3-7\alpha_4=-1.
\end{cases}
 \]
Using the second equation to eliminate $\alpha_1$ from the first and the third equation, we obtain 
 \[
\begin{cases}
3\alpha_2+6\alpha_3-3\alpha_4=-2,\\
10\alpha_2+20\alpha_3-10\alpha_4=-4,
\end{cases}
 \]
which has no solutions. Hence, no linear combination of the form \eqref{eq:linear-combination} exist, and the coset of $w$ is a nonzero element in $\WLie(x,y)/(r,s)$ that vanishes in any Novikov envelope of that algebra.
\end{proof}

\section{Invariant differential operators on the affine line}

As we suggested in the introduction, equivariant multilinear operations on the adjoint module of the Witt Lie algebra play a central role in this paper. However, some of our results are valid in larger generality. Namely, a central question about the $W_1$-modules of densities, going back to the work of Feigin and Fuchs \cite{MR554429} and Kirillov \cite{MR611158}, is to describe explicitly the vector spaces of all equivariant multilinear maps 
$\Hom_{W_1}(F_{\alpha_1}\otimes\cdots\otimes F_{\alpha_n},F_\beta)$.  
Elements of these vector spaces are known as invariant differential operators on the affine line. In this section, we prove several results about them.

\subsection{Affine-invariant multilinear operations on density modules}

Let us record an important observation concerning an intermediate invariance property, where we only require equivariance with respect to the Lie  subalgebra $\mathfrak{b}\subset W_1$ spanned by the vector fields $X_{-1}=\partial$ and $X_0=x\partial$; this is the Lie algebra of the Lie group of affine transformations.

\begin{proposition}\label{prop:borel-inv}
The vector space
\begin{equation}\label{eq:borel-inv}
T_{\alpha_1,\ldots,\alpha_n}^\beta:=\Hom_{\mathfrak{b}}(F_{\alpha_1}\otimes\cdots\otimes F_{\alpha_n},F_\beta)   
\end{equation}
vanishes unless $\beta-\alpha_1-\cdots-\alpha_n\in\mathbb{N}$. In the latter case, the vector space
 \[
\Hom_{\mathfrak{b}}(F_{\alpha_1}\otimes\cdots\otimes F_{\alpha_n},F_\beta)  
 \]
is isomorphic to the vector space of homogeneous polynomials in variables $t_1,\ldots,t_n$ of total degree $\beta-\alpha_1-\cdots-\alpha_n$.
Moreover, if we associate to each such polynomial $P(t_1,\ldots,t_n)=\sum c_{d_1,\ldots,d_n}t_1^{d_1}\cdots t_n^{d_n}$ the map 
\begin{equation}\label{eq:difop}
\Phi_P\colon f_1(x)dx^{\alpha_1}\otimes\cdots\otimes f_n(x)dx^{\alpha_n}\mapsto 
\sum c_{d_1,\ldots,d_n} f_1^{(d_1)}(x)\cdots f_n^{(d_n)}(x)dx^\beta,   
\end{equation}
then the composition of equivariant multilinear maps
\begin{equation}
\circ_i \colon T_{\alpha_1,\ldots,\alpha_n}^\beta\otimes T_{\gamma_1,\ldots,\gamma_p}^{\alpha_i}\to T_{\alpha_1,\ldots,\alpha_{i-1},\gamma_1,\ldots,\gamma_p,\alpha_{i+1},\ldots,\alpha_n}^\beta
\end{equation}
corresponds, on the level of polynomials, to the map
\begin{equation}\label{eq:comp-borel-inv}
(P\circ_i Q)(t_1,\ldots,t_{n+p-1}):=P(t_1,\ldots,t_{i-1},t_i+\cdots+t_{i+p-1},t_{i+p},\ldots,t_{n+p-1})Q(t_i,\ldots,t_{i+p-1}).
\end{equation}
\end{proposition}

\begin{proof}
Let us note that for every $\alpha$, we have $X_{-1}(f(x)dx^\alpha)=f'(x)dx^\alpha$, so the module $F_\alpha$ is co-free for the Lie subalgebra of $W_1$ generated by $X_{-1}$ (strictly speaking, co-free in the category of modules on which $X_{-1}$ acts nilpotently, which is sufficient for what follows), and does not depend on $\alpha$. On the other hand, the element $X_0$ multiplies each monomial $x^kdx^\alpha$ by the scalar $k+\alpha$, so the module $F_\alpha$ has a basis of eigenvectors of $X_0$. The co-freeness implies that 
 \[
\Hom_{\mathfrak{b}}(F_{\alpha_1}\otimes\cdots\otimes F_{\alpha_n},F_\beta)\cong \Hom_{\bk X_0}(F_{\alpha_1}\otimes\cdots\otimes F_{\alpha_n},dx^\beta),
 \]
and it remains to note that an $X_0$-equivariant linear map does not vanish on the vector 
 \[
x^{d_1}dx^{\alpha_1}\otimes\cdots\otimes x^{d_n}dx^{\alpha_n} \in F_{\alpha_1}\otimes\cdots\otimes F_{\alpha_n}
 \]
only if $d_1+\alpha_1+\cdots+d_n+\alpha_n=\beta$, that is,  $\beta-\alpha_1-\cdots-\alpha_n=d_1+\cdots+d_n\in\mathbb{N}$. This proves the first assertion: our vector space has a basis of ``multi-indices'' $(d_1,\ldots,d_n)$ with $d_1+\cdots+d_n=\beta-\alpha_1-\cdots-\alpha_n$, which identifies it with the vector space of homogeneous polynomials of the corresponding total degree. 

To prove the second assertion, we note that if $\mathbf{m}:=t_1^{d_1}\cdots t_n^{d_n}$ and $\mathbf{m}':=t_1^{e_1}\cdots t_p^{e_p}$ are two monomials, then $\Phi_{\mathbf{m}}\circ_i \Phi_{\mathbf{m}'}$ is the composite map
\begin{multline*}
f_1(x)dx^{\alpha_1}\otimes\cdots\cdots f_{i-1}(x)dx^{\alpha_{i-1}}\otimes g_1(x)dx^{\gamma_1}\otimes\cdots\otimes g_p(x)dx^{\gamma_p}\otimes f_{i+1}(x)dx^{\alpha_{i+1}} \cdots \otimes f_{n}(x)dx^{\alpha_{n}}\\
\mapsto 
f_1(x)dx^{\alpha_1}\otimes\cdots\cdots f_{i-1}(x)dx^{\alpha_{i-1}}\otimes g_1^{(e_1)}(x)\cdots g_p^{(e_p)}(x)dx^{\alpha_i}\otimes f_{i+1}(x)dx^{\alpha_{i+1}} \cdots \otimes f_{n}(x)dx^{\alpha_{n}}\\
\mapsto 
f_1^{(d_1)}(x)\cdots f_{i-1}^{(d_{i-1})}(x) (g_1^{(e_1)}(x)\cdots g_p^{(e_p)}(x))^{(d_i)} f_{i+1}^{(d_{i+1})}(x) f_{n}^{(d_n)}(x)dx^{\beta},
\end{multline*}
and the statement follows by induction on $d_i$ from the formula
 \[
(g_1^{(e_1)}(x)\cdots g_p^{(e_p)}(x))'=\sum_{s=1}^p g_1^{(e_1)}(x)\cdots g_s^{(e_s+1)}(x)\cdots g_p^{(e_p)}(x),    
 \]
corresponding, on the level of polynomials, to the multiplication by $t_i+\cdots+t_{i+p-1}$.
\end{proof}

We shall refer to the correspondence between homogeneous polynomials and equivariant multilinear maps recorded in Proposition \ref{prop:borel-inv} as the \emph{symbol} calculus for $\mathfrak{b}$-equivariant multilinear maps: the polynomial $P$ is the symbol of the multilinear map $\Phi_P$.\\

This result has an important corollary. 
\begin{corollary}\label{cor:b-Nov}
The operad $\End^{\mathfrak{b}}_{W_1}$ is isomorphic to the operad $\Nov$ of Novikov algebras.
\end{corollary}

\begin{proof}
Recall that the adjoint $W_1$-module is precisely the density module $F_{-1}$. Thus, 
 \[
\End^{\mathfrak{b}}_{W_1}(n)=\Hom_{\mathfrak{b}}(W_1^{\otimes n},W_1)= \Hom_{\mathfrak{b}}(F_{-1}^{\otimes n},F_{-1}),   
 \]
which we identified with the vector space of polynomials of degree $-1-(-1+\cdots-1)=n-1$. Moreover, the explicit computation of substitution maps $\circ_i$ in Proposition \ref{prop:borel-inv} identifies this space with differential monomials of differential order $n-1$ with the correct operad structure.  
\end{proof}

It is interesting to compare this result with work of McLachlan, Modin, Munthe-Kaas, and Verdier~\cite{MR3510021} that can be reformulated by saying that the operad of affine-equivariant multilinear maps on the adjoint module of $W_n$ for \emph{unspecified $n$} is isomorphic to the operad of pre-Lie algebras, defined by first of the two identities of Novikov algebras alone. There does not seem to be a similar ``nice'' explicit characterization of the operad of affine-equivariant multilinear maps on the adjoint module of $W_n$ for a fixed $n\ne 1$; we refer the reader to work of Bruned and Laubie \cite{BL2025} for a discussion of this matter.

A natural further question would be to consider another intermediate invariance property, where we only require equivariance with respect to the Lie  subalgebra $\mathfrak{sl}_2$ of the Witt Lie algebra spanned by the vector fields $X_{-1}=\partial$, $X_0=x\partial$, and $X_1=x^2\partial$; this is the Lie algebra of the Lie group of Möbius transformation of the projective line. It turns out that the corresponding operad is a rather non-trivial suboperad of the operad of Novikov algebras generated by the Lie bracket and one extra ternary operation; the proof of this result will be presented in a forthcoming paper.

\subsection{Three characterizations of invariant differential operators}\label{sec:3charact}

In this section, we use the result of Proposition \ref{prop:borel-inv} to give three equivalent characterizations of $W_1$-equivariant multilinear operations on density modules. To state it, let us introduce an important bit of notation. For $\overline\alpha=(\alpha_1,\ldots,\alpha_n)\in\bk^n$ and a positive integer $p$, we define a differential operator $\cD_p(\overline\alpha)$ on the ring of polynomials in variables $t_1,\ldots,t_n$ by the formula
 \[
\cD_p(\overline\alpha)=\sum_{i=1}^n\left(t_i\partial_i^{p+1}+(p+1)\alpha_i\partial_i^p\right) .   
 \]
We also denote by $\Delta_i(P)=P(t_1,\ldots,t_i+1,\ldots,t_n)-P(t_1,\ldots,t_n)$ the finite difference operator with respect to the variable $t_i$.

\begin{proposition}\label{prop:3conds}
Under the identification of $\Hom_{\mathfrak{b}}(F_{\alpha_1}\otimes\cdots\otimes F_{\alpha_n},F_{\alpha_1+\cdots+\alpha_n+k})$ with the vector space of homogeneous polynomials in $n$ variables $t_1$, \ldots, $t_n$ of degree $k$, the polynomials $P$ that are symbols of $W_1$-equivariant  maps are precisely those characterized by the following equivalent conditions:
\begin{enumerate}
\item the differential condition:
\begin{equation}\label{eq:invDifOp-PDE}
\cD_p(\overline\alpha)(P)=0, \quad p>0,
\end{equation}
\item the difference condition: \begin{equation}\label{eq:invDifOp-difference}
\sum_{i=1}^n(t_i+\alpha_i)\Delta_i(P)=k P,
\end{equation}
\item the vanishing condition: $P$ vanishes at all shifted integer points of the simplex 
\begin{equation}\label{eq:simplex}
V_k(\alpha_1,\ldots,\alpha_n):=
\left\{
(-\alpha_1-s_1,\dots,-\alpha_n-s_n)\in\bk^n\colon 
s_i\in\mathbb{N},\quad s_1+\cdots+s_n\le k-1
\right\}.
\end{equation}
\end{enumerate}
\end{proposition}
 
\begin{proof}
The differential condition is well known, see \cite[Lemma 4]{MR611158} or \cite[Sec.~2.2]{MR659165}; and it can also be easily extracted from the proof of Proposition \ref{prop:borel-inv} above. Indeed, $X_0$-equivariance means that all $W_1$-equivariant multilinear maps are contained in the graded $\Hom$-space, which allows us to write
\begin{gather*}
\Hom_{\mathfrak{b}}(F_{\alpha_1}\otimes\cdots\otimes F_{\alpha_n},F_\beta)\cong\Hom_{\mathfrak{b}}(G_\beta,G_{\alpha_1}\otimes\cdots\otimes G_{\alpha_n}),\\
\Hom_{W_1}(F_{\alpha_1}\otimes\cdots\otimes F_{\alpha_n},F_\beta)\cong\Hom_{W_1}(G_\beta,G_{\alpha_1}\otimes\cdots\otimes G_{\alpha_n}),     
\end{gather*}
where $G_\alpha=(F_\alpha)^\vee$ is the graded dual of the module $F_\alpha$. In the proof of Proposition \ref{prop:borel-inv}, we used the fact that the module $F_\alpha$ is co-free for the subalgebra spanned by $X_{-1}$; the module $G_\alpha$ is free. Moreover, the entire Lie algebra $W_1$ acts on that module by the differential operators 
 \[
X_i^\vee:=-t\partial^{i+1}-(i+1)\alpha\partial^i, \quad i\ge -1. 
 \] 
In particular, $X_{-1}^\vee$ acts by multiplication by $-t$, which immediately shows that as a $\mathfrak{b}$-module, $G_\alpha$ is induced from the one-dimensional $\bk X_0^\vee$-module $\bk v_\alpha$ where the action is defined by $X_0^\vee(v_\alpha)=-\alpha v_\alpha$.  
Furthermore,
 \[
\Hom_{\mathfrak{b}}(G_\beta,G_{\alpha_1}\otimes\cdots\otimes G_{\alpha_n})\cong 
\Hom_{\bk X_0^\vee}(\bk v_\beta,G_{\alpha_1}\otimes\cdots\otimes G_{\alpha_n}), 
 \]
and, once we take into account the $\mathfrak{b}$-equivariance condition and replace $\mathfrak{b}$-equivariant maps by their symbols, the remaining part of $W_1$-equivariance manifestly asserts that a symbol of an equivariant operator is annihilated by all operators $\cD_p(\overline\alpha)$.

The difference condition is proved similarly to \cite[Lemma 4']{MR611158} as follows. Using the Taylor formula
 \[
\Delta_i(P)=\sum_{p\ge 1}\frac{\partial_i^p(P)}{p!}, 
 \]
we note that    
\begin{multline*}
\sum_{i=1}^n(t_i+\alpha_i)\Delta_i(P)=\sum_{i=1}^n (t_i+\alpha_i)\sum_{p\ge 1}\frac{\partial_i^p(P)}{p!}\\
=\sum_{i=1}^n t_i\partial_i(P)+\sum_{p\ge 1}\frac{1}{(p+1)!}\left(\sum_{i=1}^n
t_i\partial_i^{p+1}+(p+1)\alpha_i\partial_i^p\right)(P)
=k P+\sum_{p\ge 1}\frac{1}{(p+1)!}\cD_p(\overline\alpha)(P),
\end{multline*}
so the difference condition is equivalent to 
 \[
\sum_{p\ge 1} \frac{1}{(p+1)!}\cD_p(\overline\alpha)(P)=0, 
 \]
which, because the polynomial $P$ is assumed homogeneous and because each operator $\cD_p(\overline\alpha)$ decreases degree by $p$, implies that $\cD_p(\overline\alpha)(P)=0$ for all $p\ge 1$, which is precisely the differential condition. 

Let us show that the difference condition implies the vanishing condition. As a first step, we remark that for $n=1$ the rising factorial polynomial 
 \[
f_{k_1,\alpha_1}(t_1)=(t_1+\alpha_1)(t_1+\alpha_1+1)\cdots (t_1+\alpha_1+k_1-1) 
 \]
is an eigenvector of $(t_1+\alpha_1)\Delta_1$ with the eigenvalue $k_1$, since $(t_1+\alpha_1)\Delta_1(f_{k_1,\alpha_1}(t_1))$ is given by
 \[ 
(t_1+\alpha_1)(t_1+\alpha_1+1)\cdots (t_1+\alpha_1+k_1)-(t_1+\alpha_1)(t_1+\alpha_1)(t_1+\alpha_1+1)\cdots (t_1+\alpha_1+k_1-1)=
k_1f_{k_1,\alpha_1}(t_1).
 \]
It follows that for $\mathbf{k}=(k_1,\ldots,k_n)\in\mathbb{N}^n$, the polynomial
\begin{equation}
f_{\mathbf{k},\overline\alpha}(t_1,\ldots,t_n):=
f_{k_1,\alpha_1}(t_1)\cdots f_{k_n,\alpha_n}(t_n) 
\end{equation}
is an eigenvector of the operator
 \[
E^\Delta(\overline\alpha):=\sum_{i=1}^n(t_i+\alpha_i)\Delta_i 
 \]
with the eigenvalue $k_1+\cdots+k_n$. Since the polynomials 
 \[
\left\{f_{\mathbf{k},\overline\alpha}(t_1,\ldots,t_n)\right\}_{\mathbf{k}\in\mathbb{N}^n}
 \]
are related to the standard monomial basis of the polynomial ring by an unitriangular matrix, they form a basis of $\bk[t_1,\ldots,t_n]$, and we have just observed that the operator $E^\Delta(\overline\alpha)$ is diagonal in that basis. Thus, Condition \eqref{eq:invDifOp-difference} is equivalent to the fact that the polynomial $P$ is a linear combination of polynomials $f_{\mathbf{k},\overline\alpha}(t_1,\ldots,t_n)$ with $k_1+\cdots+k_n=k$. However, if $\mathbf{k}\in\mathbb{N}^n$ satisfies $k_1+\cdots+k_n=k$, then  $f_{\mathbf{k},\overline\alpha}(t_1,\ldots,t_n)$ vanishes on $V_k(\overline\alpha)$. Indeed, 
 \[
f_{\mathbf{k},\overline\alpha}(-\alpha_1-s_1,\ldots,-\alpha_n-s_n)= \prod_{i=1}^n\prod_{j=0}^{k_i-1}(-s_i+j).
 \]
If $s_1+\cdots+s_n\le k-1=k_1+\cdots+k_n-1$, then for some $i$ one has
$0\le s_i<k_i$, and the factor above corresponding to $j=s_i$ vanishes. 

Conversely, suppose that the vanishing condition holds for a homogeneous polynomial $P$ of degree $k$. We note that this immediately implies that for each $i\in\{1,\ldots,n\}$ the polynomial 
 \[
(t_i+\alpha_i)P(t_1,\ldots,t_i+1,\ldots,t_n) 
 \]
also vanishes on $V_k(\alpha_1,\ldots,\alpha_n)$: indeed, if 
 \[
(-\alpha_1-s_1,\ldots,-\alpha_n-s_n)\in V_k(\alpha_1,\ldots,\alpha_n), 
 \]
then for $s_i=0$ already $t_i+\alpha_i$ vanishes, and for $s_i>0$, we note that the evaluation of $P(t_1,\ldots,t_i+1,\ldots,t_n)$ at the point $(-\alpha_1-s_1,\ldots,-\alpha_n-s_n)$ is the same as the evaluation of $P$ at the point 
$(-\alpha_1-s_1,\ldots,-\alpha_i-(s_i-1),\ldots,-\alpha_n-s_n)$, which manifestly belong to $V_k(\alpha_1,\ldots,\alpha_n)$.
Thus the polynomial
 \[
\sum_{p=1}^n(t_i+\alpha_i)\Delta_i(P)- kP.  
 \]
vanishes on $V_k(\alpha_1,\ldots,\alpha_n)$ and has degree at most $k-1$. It remains to use the elementary interpolation fact that a polynomial
of total degree at most $d$ which vanishes at all points of $\mathbb{N}^n$ with sum of coordinates not exceeding $d$ is identically zero.
\end{proof} 

Let us record a particular case of our result for $\alpha_1=\cdots=\alpha_n=\beta=-1$ that will be crucial for the main result of the paper.

\begin{corollary}\label{cor:inv-vanishing}
Under the identification of $\Hom_{\mathfrak{b}}(W_1^{\otimes n},W_1)$ with the vector space of homogeneous polynomials in $n$ variables $t_1$, \ldots, $t_n$ of degree $n-1$, the polynomials that are symbols of multilinear $W_1$-equivariant maps are precisely those that vanish at all points of the shifted simplex 
\begin{equation}\label{eq:simplex-ad}
V_{n-1}:=
\left\{
(1-s_1,\dots,1-s_n)\in\mathbb{Z}^n\colon 
s_i\in\mathbb{N},\quad s_1+\cdots+s_n\le n-2
\right\}.
\end{equation}
\end{corollary}

As we were informed by Molev, the vanishing conditions that appear in this and the following sections were first discovered by Kontsevich in the context of identities of Lie algebras of vector fields and recorded in the paper \cite{MR713969} by Kirillov and Kontsevich; to the best of our knowledge, their meaning in the context of invariant differential operators was never realized. 

\section{Wronskian Lie polynomials and equivariant multilinear maps}

We shall now use the beautiful results of Molev \cite{Mo} about the Lie suboperad of $\End_{W_1}$ to prove the main result of this paper. In fact, a large part of the below proof is a careful adaptation of the ingenious proof of \cite[Th.~1]{Mo}; moreover, our work was easier since we had to deal with radical ideals, encountering the relevant vanishing conditions directly, and not through a concrete set of polynomials satisfying those vanishing conditions. Nevertheless, while the argument is similar, it is not the same, and we chose to spell it out in detail; additionally, we work with vector fields with polynomial coefficients which made certain arguments of \cite{Mo} unavailable. 

\begin{theorem}\label{th:LieNovInv}
The operad of Wronskian Lie algebras is isomorphic to the operad of equivariant multilinear operations on the adjoint module of the Witt Lie algebra. 
\end{theorem}

\begin{proof}
Let us note that the Jacobi identity, written in the form $[a,[b,c]]=[[a,b],c]+[b,[a,c]]$,    
implies that the Lie bracket on any Lie algebra $\mathfrak{g}$ is a $\mathfrak{g}$-equivariant map from $\mathfrak{g}^{\otimes 2}$ to $\mathfrak{g}$. Thus the image of the canonical map of operads $\Lie\to\End_{W_1}$ is contained in the suboperad $\End_{W_1}^{W_1}$ of $W_1$-equivariant multilinear maps. By Proposition \ref{prop:suboperads}, that image is isomorphic to the operad $\WLie$, and from now on we shall identify it with $\WLie$. It follows that in order to prove the assertion of the theorem, it is enough to show that for all $n\ge 1$, the dimension of $\WLie(n)$ is equal to the dimension of $\End_{W_1}^{W_1}(n)$. In fact, it will be easier to prove that those two vector spaces are isomorphic to $\bk S_n$-modules.

As we saw in the proof of Corollary \ref{cor:Molev-WLie}, the $S_n$-module structure of $\WLie(n)$ can be extracted from results of Molev \cite{Mo}. We shall establish that the $S_n$-module $\End_{W_1}^{W_1}(n)$ has the same structure. 

\textbf{Step 1: a representation-theoretic reduction. } 
Rather than compare the two $S_n$-modules directly, we shall examine their elements that transform according to particular symmetry types after restriction to suitable subgroups of $S_n$. Let us make this informal idea precise. Let $\lambda\vdash n$ be a partition of $n$, which we represent by a Young diagram. Let $l$ be the number of columns of height at least three in $\lambda$, which we shall refer to as ``tall'' columns; we denote their heights $c_1\ge \ldots\ge c_l\ge 3$. The remaining part of the diagram consists of two rows of lengths $a\ge b\ge 0$. Let us number the boxes of $\lambda$ by $1,\ldots, n$ as follows: the boxes in the first row outside of the tall columns numbered $1,\ldots,a$, the boxes in the second row outside of the tall columns numbered $a+1,\ldots,a+b$, the boxes in the first tall columns $a+b+1,\ldots,a+b+c_1$, etc. This numbering defines a subgroup $H_\lambda$ of $S_n$ isomorphic $S_a\times S_b\times S_{c_1}\times \cdots S_{c_l}$ which permutes entries each the first two rows outside the tall columns and within each of the tall columns. For an $S_n$-module $V$, we may consider the vector space $V(a,b;c_1,\ldots,c_l)$ consisting of all elements that are invariant with respect to $S_a\times S_b$ and transform according to the sign representation with respect to $S_{c_1}\times \cdots S_{c_l}$. It turns out that if $V(a,b;c_1,\ldots,c_l)$ is equal to zero for every partition $\lambda\vdash n$, then $V=0$.
A stronger result is established in \cite[Lemma 1]{Mo}, but we can also give a short self-contained proof: $V(a,b;c_1,\ldots,c_l)$ is precisely the isotypic component of the $S_a\times S_b\times S_{c_1}\times \cdots S_{c_l}$-module
 \[
U_\lambda:=\mathrm{triv}\boxtimes\mathrm{triv}\boxtimes\mathrm{sgn}\boxtimes\cdots\boxtimes\mathrm{sgn}     
 \]  
in $V$, and by the Frobenius reciprocity we have
 \[
\dim\Hom_{H_\lambda}(U_\lambda,\Res_{H_\lambda}^{S_n}V)= 
\dim\Hom_{H_\lambda}(\Ind_{H_\lambda}^{S_n}U_\lambda,V). 
 \]
Under the characteristic map producing symmetric functions from $S_n$-modules, \cite[Sec.~I.7]{MR3443860}, the $S_n$-module $\Ind_{H_\lambda}^{S_n}$ corresponds to the symmetric function $h_ah_be_{c_1}\cdots e_{c_l}$, and the Pieri rule \cite[Sec.~I.5]{MR3443860} shows that the Schur function expansion of this symmetric function contains $s_\lambda$ with a strictly positive coefficient. Thus, if $V$ were a non-zero module, we would have that $V(a,b;c_1,\ldots,c_l)\ne 0$ for some $\lambda$. 

Recall that by Corollary \ref{cor:inv-vanishing}, $\End_{W_1}^{W_1}(n)$ may be identified with the $S_n$-module $\calV(n)$ of all polynomials of degree $n-1$ that vanish on all points of the shifted simplex 
 \[
V_{n-1}=\left\{
(1-s_1,\ldots,1-s_n)\in\mathbb{Z}^n\colon 
s_i\in\mathbb{N},\quad s_1+\cdots+s_n\le n-2
\right\}.
 \]
The category of $\bk S_n$-modules is semisimple and $\WLie(n)\subseteq\End_{W_1}^{W_1}(n)\cong \calV(n)$, therefore it is enough to prove that the $\bk S_n$-module $\calM(n):=\calV(n)/\WLie(n)$ satisfies
$\calM(n)(a,b;c_1,\ldots,c_l)=0$ for all $\lambda\vdash n$, and that would immediately follow from the inequalities
\begin{equation}\label{eq:ineq-needed}
\dim\calV(n)(a,b;c_1,\ldots,c_l)\le \dim\WLie(n)(a,b;c_1,\ldots,c_l) 
\end{equation}
for all partitions $\lambda\vdash n$. Establishing these inequalities is precisely what we shall do now. 

\textbf{Step 2: the mixed symmetry components of $\WLie(n)$. }
Let us deduce from Corollary \ref{cor:Molev-WLie} a precise formula for the right hand side of \eqref{eq:ineq-needed}. First, we note that the isotypic component of the $S_a\times S_b\times S_{c_1}\times \cdots S_{c_l}$-module $U_\lambda$ in $\bk[t_1,\ldots,t_n]$ consists of all polynomials that are symmetric in $t_1,\ldots,t_a$, and separately in $t_{a+1},\ldots,t_{a+b}$, and skew-symmetric in the groups of variables of each tall column of $\lambda$. Since the vector space of skew-symmetric polynomials in several variables is, as a module over the ring of symmetric polynomial, a free rank one module generated by the product of pairwise differences of variables, we see that the Hilbert series of the isotypic component in question is given by the formula
 \[
H_a(q)H_b(q)\prod_{i=1}^lq^{\binom{c_i}2}H_{c_i}(q), 
 \]
where $H_u(q)=((1-q)(1-q^2)\cdots(1-q^u))^{-1}$     
is the Hilbert series of the ring of symmetric functions in $u$ variables. By a classical result of Chevalley~\cite{MR72877}, the ring $\bk[t_1,\ldots,t_n]$ is a free module over the subring $\bk[t_1,\ldots,t_n]^{S_n}$ of $S_n$-invariants. Clearly, the space of generators of that module is naturally identified with the coinvariant ring $\calH_n$, and we an isomorphism of $S_n$-modules 
$\bk[t_1,\ldots,t_n]\cong \bk[t_1,\ldots,t_n]^{S_n}\otimes  \calH_n$. Since $\bk[t_1,\ldots,t_n]^{S_n}$ is a trivial graded $S_n$-module with the Hilbert series $H_n(q)$, it follows that the isotypic component of the $S_a\times S_b\times S_{c_1}\times \cdots S_{c_l}$-module $U_\lambda$ in $\calH_n$ has the Hilbert series 
 \[
\frac{1}{H_n(q)}H_a(q)H_b(q)\prod_{i=1}^lq^{\binom{c_i}2}H_{c_i}(q)=\prod_{i=1}^n(1-q^i)H_a(q)H_b(q)\prod_{i=1}^lq^{\binom{c_i}2}H_{c_i}(q). 
 \]
Thus, it is an immediate consequence of Corollary \ref{cor:Molev-WLie} that
 \[
\dim\WLie(n)(a,b;c_1,\ldots,c_l)=[q^{n-1}]\prod_{i=1}^n(1-q^i)H_a(q)H_b(q)\prod_{i=1}^lq^{\binom{c_i}2}H_{c_i}(q) . 
 \]
This means that our result will be proved once we show that 
\begin{equation}\label{eq:transformed-ineq}
\dim\calV(n)(a,b;c_1,\ldots,c_l)\le [q^{n-1}]\prod_{i=1}^n(1-q^i)H_a(q)H_b(q)\prod_{i=1}^lq^{\binom{c_i}2}H_{c_i}(q).    
\end{equation}
Let us remark that \eqref{eq:transformed-ineq} is straightforward for Young diagrams $\lambda$ for which all columns are of height at least three. Indeed, in this case,  each element from $\calV(n)$ is skew-symmetric in the groups of $c_1,\ldots,c_l$ variables with $c_1+\cdots+c_l=n$. This implies that it is divisible by products of differences of variables in each group, and so its degree is at least 
 \[
\sum_{i=1}^l\binom{c_i}2\ge \sum_{i=1}^n c_i=n, 
 \]
since each $c_i$ is at least three. This latter inequality shows that both sides of \eqref{eq:transformed-ineq} are equal to zero. 

\textbf{Step 3: getting rid of the $a$-block of variables. }
Consider now a Young diagram $\lambda$ which has some columns of height at most two; in the notation introduced above, this means that $a\ge 1$. Let us show that if we specialize the first $a$ variables to zero, no information is lost.  After this specialization, only the second symmetric block and the skew-symmetric tall-column blocks remain, and their possible leading terms will be estimated by an interpolation argument.

The following lemma is an algebraic version of the elegant idea of Kirillov, Kontsevich and Molev \cite[Part~II]{MR753873} that when working with identities satisfied by generic vector fields on the affine line, one may assume that one of those vector fields is simply $\partial$. In the algebraic category, doing an analytic straightening of a vector field is not an option, and our argument below is precisely an algebraic counterpart of such analytic argument.

\begin{lemma}\label{lm:evaluation-a}
Setting $t_1=\cdots=t_a=0$ defines an injective map 
$\calV(n)(a,b;c_1,\ldots,c_l)\hookrightarrow\bk[t_{a+1},\ldots,t_n]$.   
\end{lemma}

\begin{proof}
Let $P_\Psi=\sum_{i_1+\cdots+i_n=n-1}c_{i_1,\ldots,i_n}t_1^{i_1}\cdots t_n^{i_n}$, which, as the notation suggests, we view as the symbol of the differential polynomial 
 \[
\Psi=\sum_{i_1+\cdots+i_n=n-1}c_{i_1,\ldots,i_n}u_1^{(i_1)}\cdots u_n^{(i_n)}. 
 \]
Because of our assumption on the degree, $\Psi$ is an element of the ring of differential polynomials in $u_1$, \ldots, $u_n$ of order at most $n-1$ in each variable. Suppose that the evaluation of the polynomial $P$ at $t_1=\cdots=t_a=0$ is equal to zero. That evaluation corresponds to setting $u_1=\cdots=u_a=u$ in $\Psi$, retaining the differential monomials not containing $u^{(p)}$ with $p>0$, and deleting $u^{(0)}$. Now let us briefly view the ring of differential polynomials in $u_1$, \ldots, $u_n$ of order at most $n-1$ in each variable as simply the ring of polynomials 
 \[
\bk[u_i^{(j)}\colon 1\le i\le n, 0\le j\le n-1].   
 \]
Let us consider the derivations $\delta_k$ of this ring defined by 
 \[
\delta_k(u_i^{(j)}=
\begin{cases}
\left(\binom{j}{k+1}-\binom{j}{k}\right)u_i^{(j-k)}, j\ge k\\
\qquad\qquad 0, \qquad\qquad j<k,
\end{cases}
 \]
In fact, since
 \[
\frac{1}{(k+1)!}(t\partial^{k+1}-(k+1)\partial^k)(t^j)=
\left(\binom{j}{k+1}-\binom{j}{k}\right)t^{j-k},
 \]
these derivations correspond to the action of the operators from Section \ref{sec:3charact}:
 \[
P_{\delta_k(\Psi)}=\frac{1}{(k+1)!}\cD_k(\overline\alpha)(P_\Psi) 
 \]
for $\alpha=(-1,-1,\ldots,-1)$. However, unlike the differential operators $\cD_k(\overline\alpha)$ acting on polynomials $t_i$, in the chosen context of differential monomials our operators $\delta_k$ are locally nilpotent derivations, and hence can be exponentiated to automorphisms $\exp(z\delta_k)$ of the ring $\bk[u_i^{(j)}\colon 1\le i\le n, 0\le j\le n-1]$. We note that if $P_\Psi$ vanishes on $V_{n-1}$, then $\cD_k(\overline\alpha)(P_\Psi)=0$ for all $k\ge 1$, meaning that $\delta_k(\Psi)=0$ for all $k$, 
and therefore $\Psi$ is invariant under all automorphisms $\exp(z\delta_k)$. 

Let us now set $u_1^{(j)}=\cdots=u_a^{(j)}=u^{(j)}$ in $\Psi$ for each $j=0,\ldots,n-1$; we denote by $\Psi^=$ the resulting polynomial, which is now an element of 
 \[
\bk[u^{(j)}, u_i^{(j)}\colon a+1\le i\le n, 0\le j\le n-1]. 
 \]
Note that the passage from $\Psi$ to $\Psi^=$ commutes with the action of each $\delta_k$, since the action of these derivations on the variable $u_i^{(j)}$ does not depend on $i$. Hence $\delta_k(\Psi^=)=0$ and $\Psi^=$ is invariant under all automorphisms $\exp(z\delta_k)$. 

Let us now completely ignore the differential ring structure and view $u^{(j)}, u_i^{(j)}$ as formal variables which are coordinates on an affine space of dimension $n(n-a+1)$. Let us take a point $\overline{u}$ of that space where $\overline{u}^{(0)}\ne 0$. Since 
 \[
\delta_k(u^{(k)}=-u^{(0)}, \qquad \delta_k(u^{(r)})=0 \text{ for } r<k,  
 \]
applying the automorphism $\exp(z\delta_k)$ with $z=\frac{u^{(k)}}{u^{(0)}}$ moves from the point $\overline{u}$ to the point $\overline{v}$ for which $\overline{v}^{(k)}=0$ and $\overline{v}^{(r)}=\overline{v}^{(r)}$ for $r<k$. Let us apply that procedure consecutively for $k=1$, then for $k=2$, etc., up to $k=n-1$. As a result, we shall move from point $\overline{u}$ to some point $\overline{\overline{u}}$ for which $\overline{\overline{u}}^{(r)}=0$ for $1\le r\le n-1$. Since $\Psi^=$ is of degree $a$ in the variables $u^{(s)}$, its value at the point $\overline{\overline{u}}$ is equal to the product of $(\overline{\overline{u}}^{(0)})^a$ and a polynomial in all the remaining variables, and that latter polynomial is precisely the symbol of the evaluation of $P_\Psi$ at $t_1=\cdots=t_a=0$. By our assumption, that evaluation is zero, and since $\Psi^=$ is invariant under the automorphisms we used, we conclude that $\Psi^=$ vanishes at all points $\overline{u}$ where $\overline{u}^{(0)}\ne 0$, and therefore vanishes identically. Finally, since $\Psi$ is symmetric in its first $a$ arguments, it can be recovered from $\Psi^=$ by polarization. 
\end{proof}

\textbf{Step 4: estimating the dimension after specialization. }
We have an injective map
 \[
\calV(n)(a,b;c_1,\ldots,c_l)\to\bk[t_{a+1},\ldots,t_n]\cong \bk[t_{a+1},\ldots,t_{a+b}]\otimes \bk[t_{a+b+1},\ldots,t_n]   
 \]
obtained by setting $t_1=\cdots=t_a=0$.  In fact, we have an obvious constraint on the image of that map. By definition, a polynomial $P\in \calV(n)$ vanishes at all points of the shifted simplex
 \[
V_{n-1}=\left\{
(1-s_1,\ldots,1-s_n)\in\mathbb{Z}^n\colon 
s_i\in\mathbb{N},\quad s_1+\cdots+s_n\le n-2
\right\}. 
 \]
We note that each point of that simplex of the form
 \[
(0,\ldots,0,1-s_{a+1},\ldots,1-s_n)=(1-1,\ldots,1-1,1-s_{a+1},\ldots,1-s_n)
 \]
gives rise to a point $(1-s_{a+1},\ldots,1-s_n)$ in the simplex of the same type in lower dimension: we have 
 \[
s_{a+1}+\cdots+s_n\le n-a-2,     
 \]
and every integer point of this lower-dimensional simplex can of course be obtained this way. Thus, the polynomial $P_0:=P|_{t_1=\cdots=t_a=0}$ vanishes at all points of that simplex. Moreover, $P\in \calV(n)(a,b;c_1,\ldots,c_l)$, implying that the polynomial $P_0$ is symmetric in variables $t_{a+1},\ldots,t_{a+b}$ and skew-symmetric in each of the groups of elements corresponding to tall columns of $\lambda$. Thus, it can be written as 
 \[
P_0=\sum_{j=1}^r B_j(t_{a+1},\ldots,t_{a+b}) \prod_{\substack{i<j,\\i,j \text{in the same tall column}}}(t_i-t_j) C_j(t_{a+b+1},\ldots,t_n),
 \]
where $B_j(t_{a+1},\ldots,t_{a+b})$ is symmetric, and $C_j(t_{a+b+1},\ldots,t_n)$ is symmetric in each of the groups of elements corresponding to tall columns of $\lambda$, all these polynomials are homogeneous, and all polynomials $C_j$ are linearly independent among themselves. 

Let us consider a partial order on monomials with respect to the total degree in $t_{a+b+1}$,\ldots, $t_n$; by the ``leading term'' of a polynomial with respect to that partial order, we shall mean the sum of all terms where the degree in $t_{a+b+1}$,\ldots, $t_n$ is the largest possible. By doing Gaussian elimination with respect to this partial order, we can find a basis of the vector space of possible polynomials $P_0$ with linearly independent leading terms. 

For a given polynomial $P_0$, let us denote by $d_{tc}$ the degree in the ``tall column variables'' $t_{a+b+1}$,\ldots, $t_n$ of its leading term. Because of the products of pairwise differences of variables from tall columns, we have
 \[
d_{tc}\ge \sum_{i=1}^l\binom{c_i}2\ge \sum_{i=1}^l c_i=n-a-b,
 \]
implying that 
 \[
m:=n-a-d_{tc}\le b.     
 \]
Following \cite{MR753873,Mo}, let us denote by $E_{b,m}$ a slightly more general shifted simplex
 \[
\left\{(1-s_{a+1},\ldots,1-s_{a+b})\in\mathbb{Z}^b \colon s_i\in\mathbb{N}, s_{a+1}+\cdots+s_{a+b}\le m-2\right\} 
 \]
(the shifted simplex $V_{n-1}$ from \eqref{eq:simplex-ad} is $E_{n,n}$). 
For any $s_{a+1},\ldots,s_{a+b}\in\mathbb{N}$ with $s_{a+1}+\cdots+s_{a+b}\le m-2$ and any $s_{a+b+1}, \ldots,s_n\in\mathbb{N}$ with $s_{a+b+1}+\cdots+s_n\le d_{tc}$, we have
 \[
s_{a+1}+\cdots+s_{a+b}+s_{a+b+1}+\cdots+s_n\le m-2+d_{tc}=n-a-2,  
 \]
so from the previous paragraph it follows that $P_0$ vanishes at the point $(1-s_{a+1},\ldots,1-s_n)$. If we substitute specific values of $s_{a+1},\ldots,s_{a+b}$ into the polynomial $P_0$ and allow $s_{a+b+1}$, \ldots, $s_n$ to vary, we obtain a polynomial in $s_{a+b+1}$, \ldots, $s_n$ of degree at most $d_{tc}$ that vanishes at all integer points of the shifted simplex 
 \[
\left\{(1-s_{a+b+1},\ldots,1-s_n)\in\mathbb{Z}^{n-a-b}\colon s_i\in\mathbb{N},  s_{a+b+1}+\cdots+s_n\le d_{tc}\right\},     
 \]
and the same interpolation argument as in the end of the proof of Proposition \ref{prop:3conds} shows that this polynomial in $s_{a+b+1}$, \ldots, $s_n$ vanishes identically. Since we assumed the polynomials $C_j$ linearly independent, this means that each polynomial $B_j$ vanishes at each point of 
$E_{b,m}$. (Strictly speaking, the preceding argument is valid for $m\ge 2$ and $n-a-b>0$ only, but for $m\le1$ this vanishing condition is empty and hence trivially true, and for $n=a+b$ we have $m=b$ and we have the only $B$-polynomial equal to $P_0$, which is known to vanish at each point of $E_{b,b}$.) 

We shall use the above argument to estimate the dimension of the vector space spanned by the polynomials $B_j$ that appear in the leading term of $P_0$, so that the degrees of the corresponding polynomials $C_j$ are all equal to $d_{tc}$. Those $B_j$ are symmetric polynomials of degree 
 \[
r=n-1-d_{tc}=a+m-1.     
 \]
Since we saw that $m\le b$, we also have $m\le a$, and therefore we have $r\ge 2m-1$.     

The following lemma is a weaker version of the statement \cite[Lemma 3]{Mo}, which is recorded in \emph{op. cit.} without a proof; it is indicated the result is similar to \cite[Th.~6]{MR753873}. Our proof is similar but is marginally simpler.

\begin{lemma}\label{lm:restriction-homogen}
For $2\le m\le b$ and $r\ge 2m-3$, every function on $E_{b,m}$ can be obtained as a restriction of a homogeneous polynomial in $x_{a+1}$, \ldots, $x_{a+b}$ of total degree $r$. 
\end{lemma}

\begin{proof}
For each point 
 \[
 \varepsilon(s)=(1-s_{a+1},\ldots,1-s_{a+b})\in E_{b,m},  
 \]
we shall construct a homogeneous polynomial of degree $2m-3$ that does not vanish at $\varepsilon(s)$ and vanishes at all other points of $E_{b,m}$. It is enough to do so for $s_{a+1}\le \cdots\le s_{a+b}$, for any other case is reduced to this by an action of a permutation. Since $m\le b$ and the defining condition of $E_{b,m}$ implies that there are at most $m-2$ non-zero elements among $s_{a+1},\ldots,s_{a+b}$, we must have $s_{a+1}=0$. Set
 \[
A_s(x)=\prod_{j=1}^{m-1}(x_{a+1}+\cdots+x_{a+j}).
 \]
We claim that $A_s(\varepsilon(s))\ne0$, whereas
 \[
A_s(\varepsilon(s'))=0 \qquad\text{if } s'_{a+1}+\cdots+s'_{a+b}\le m-2
 \text{ and }s_{a+1}'>0.
\]
Indeed, at $\varepsilon(s)$ the partial sums in the product are all positive: they initially increase by $1$, then weakly decrease, and the
last is 
 \[
m-1-s_{a+1}-\cdots-s_{a+b}\ge1.     
 \]
On the other hand, if $s'_{a+1}>0$, the first partial sum at $\varepsilon(s')$ is $1-s'_{a+1}\le0$, while the last is at least 
 \[
m-1-s_{a+1}'-\cdots-s_{a+b}'\ge1.     
 \]
These partial sums are integers, and each successive increment is at most $1$. Hence one of them must be zero.

It remains to distinguish $\varepsilon(s)$ among the points $\varepsilon(s')$ with $s'_{a+1}=0$. On this subset, $x_{a+1}=1$, and the numbers
$y_i=x_{a+1}-x_i$ for $i>{a+1}$ run through the integer simplex
 \[
\left\{(t_i)_{i>{a+1}}\colon t_i\in\mathbb N,\quad
 \sum_{i>{a+1}}t_i\le m-2\right\}.
 \]
Standard interpolation allows us to represent every function on this simplex by a (non-homogeneous) polynomial in the variables $y_i$ of total degree at most $m-2$. Homogenizing that polynomial with respect to $x_{a+1}$ therefore gives a homogeneous polynomial $B_s(x)$ of degree $m-2$ such that
 \[
 B_s(\varepsilon(s'))=\begin{cases}
 1,&s'=s,\\
 0,&s'\ne s,
 \end{cases}
 \qquad\text{whenever }s_{a+1}'=0.
 \]
Now $F_s=A_sB_s$ has degree $m-1+m-2=2m-3$ and vanishes at every point of $E_{b,m}$ other than $\varepsilon(s)$, where it is nonzero.
Thus the polynomials $F_s$ give interpolation in degree $2m-3$. Multiplying $F_s$ by $x_{a+1}^{r-(2m-3)}$ gives the same property
in every degree $r\ge 2m-3$, because $x_{a+1}(\varepsilon(s))=1$. 
\end{proof}

It follows from this lemma that every $S_b$-invariant function on $E_{b,m}$ can be obtained as a restriction of a homogeneous symmetric polynomial of degree $r$. The dimension of the vector space of $S_b$-invariant functions on $E_{b,m}$ is equal to the number of $S_b$-orbits on $E_{b,m}$, which is equal to the number of partitions of integers not exceeding $m-2$ into at most $b$ parts; for the usual notation $p_k(m)$ for the number of partitions of $m$ into at most $k$ parts, that number of orbits is equal to
$\sum_{j=0}^{m-2}p_b(j)$; here and in the following formula, we interpret a sum from zero to a negative upper bound as zero.  
Thus, the dimension of the space of homogeneous symmetric polynomials of degree $r$ that vanish on $E_{b,m}$ is equal to 
 \[
p_b(r)-\sum_{j=0}^{m-2}p_b(j)= p_b(r)-\sum_{j=0}^{r-a-1}p_b(j).
 \] 
We remark that since $r=a+m-1\le a+b-1\le 2a-1$, we have 
 \[
\prod_{j\ge a+1}(1-q^j)=1-\sum_{j=a+1}^{2a}q^j + O(q^{2a+1}),
 \]
and hence 
 \[
p_b(r)-\sum_{j=0}^{r-a-1}p_b(j)= [q^r] H_b(q)\prod_{j\ge a+1}(1-q^j).
 \]
To summarize, the total dimension of the space of polynomials $P$ in $t_1,\ldots,t_n$ of degree $n-1$ that vanish at all points of $V_{n-1}$ and have the symmetry type $(a,b;c_1,\ldots,c_l)$ is bounded from the above by the total dimension of their evaluations at $t_1=\cdots=t_a=0$, and enumerating them according to their top degree $d_{tc}$ in $t_{a+b+1}$,\ldots, $t_n$, our discussion shows that the dimension of the vector space of such polynomials does not exceed 
\begin{multline*}
[q^{n-1}] H_b(q)\prod_{j\ge a+1}(1-q^j)\prod_{i=1}^lq^{\binom{c_i}2}H_{c_i}(q)=\\
[q^{n-1}] H_b(q)\prod_{i=a+1}^n(1-q^i)\prod_{i=1}^lq^{\binom{c_i}2}H_{c_i}(q)=
[q^{n-1}] H_b(q)\prod_{i=1}^n(1-q^i)H_a(q)\prod_{i=1}^lq^{\binom{c_i}2}H_{c_i}(q),
\end{multline*}
where the first equality follows from the fact that the discrepancy between
 \[
H_b(q)\prod_{j\ge a+1}(1-q^j)\prod_{i=1}^lq^{\binom{c_i}2}H_{c_i}(q)\quad \text{ and }\quad  H_b(q)\prod_{i=a+1}^n(1-q^i)\prod_{i=1}^lq^{\binom{c_i}2}H_{c_i}(q)    
 \]
only involves factors of the form $1-q^i$ with $i>n$ which do not affect the coefficient of $q^{n-1}$, and the second equality merely uses the definition of $H_a(q)$. This completes the proof of \eqref{eq:transformed-ineq}, which, as we indicated above, is the final step in the proof of the theorem.
\end{proof}

\begin{corollary}\label{cor:criterion-symbols}
A homogeneous polynomial $P$ of degree $n-1$ is a symbol of a multilinear Lie polynomial in $\Nov(n)$ if and only if any of the following conditions is satisfied:
\begin{enumerate}
\item the differential condition:
\begin{gather}
\sum_{i=1}^n t_i\partial_i^2(P)=2\sum_{i=1}^n\partial_i(P),\label{eq:invDifOp-PDE-minus-one-1}\\
\sum_{i=1}^n t_i\partial_i^3(P)=3\sum_{i=1}^n\partial_i^2(P),\label{eq:invDifOp-PDE-minus-one-2}
\end{gather}
\item the difference condition: \begin{equation}\label{eq:invDifOp-difference-minus-one}
\sum_{i=1}^n(t_i-1)\Delta_i(P)=(n-1)P,
\end{equation}
\item the vanishing condition: $P(a)=0$ for all $a\in V_{n-1}$, where
 \[
V_{n-1}= \bigl\{(1-b_1,\ldots,1-b_n)
\mid b\in\mathbb{N}^n,\ |b|\le n-2\bigr\} .
 \]
\end{enumerate}
\end{corollary}

\begin{proof}
This is an immediate consequence of Corollary \ref{cor:b-Nov}, Proposition \ref{prop:3conds}, and Theorem \ref{th:LieNovInv}, paired with the observation that $X_1$ and $X_2$ generate the positive part of the Witt algebra, and hence it is sufficient to check the differential condition for $p=1$ and $p=2$ only.
\end{proof}

This corollary gives three finite procedures for recognizing multilinear Lie polynomials in free Novikov algebras. To apply it to a general element of a free Novikov algebra, one should multilinearize it and then convert it into its differential presentation. As a consequence, those finite procedures are rather indirect. In the following section, we shall prove a much more direct criterion for Lie polynomials in free Novikov algebras. Before we do that, we give some examples of how Corollary \ref{cor:criterion-symbols} can be applied. 

\begin{corollary}\label{cor:no-deriv}
Let $f$ and $g$ be nonzero Lie polynomials in the free Novikov algebra. 
\begin{enumerate}
\item The differential realization of $f$ is not a total derivative, i.e.\ it is not of the form $h'$ for some differential polynomial $h$.
\item The Novikov product $f\circ g$ is not a Lie polynomial.
\end{enumerate} 
\end{corollary}

\begin{proof}
To prove the first claim, we note that multilinearization preserves the properties of being a Lie polynomial and being a total derivative, so we may assume $f$ to be multilinear. It will be convenient to pass to symbols and write $f=\Psi_u$. Since applying $D$ amounts to multiplying by $t_1+\cdots+t_n$, we have $u=(t_1+\cdots+t_n)v$. The polynomial $t_1+\cdots+t_n$ does not vanish on our shifted simplex, which forces the polynomial $v$ to vanish on it. The polynomial $v(1-t_1,\ldots,1-t_n)$ is of degree $n-2$ (though no longer homogeneous) and vanishes on the standard simplex consisting of all $s\in\mathbb{N}^n$ with $|s|\le n-2$, so our usual interpolation argument gives $v=0$, meaning that $u=0$ and therefore $f=0$.

For the second claim, we note that $f\circ g=\frac12(f\circ g+g\circ f)+\frac12[f,g]$. The second term is a Lie polynomial, so $f\circ g$ is a Lie polynomial if and only if $f\circ g+g\circ f$ is a Lie polynomial. But in the differential realization, we have $f\circ g+g\circ f=(fg)'$, and the first assertion completes the proof.
\end{proof}

As established by Dzhumadildaev and Ismailov \cite[Th.~3.2]{DI2026}, Novikov elements that are total derivatives are precisely iterations of the symmetrized product $a\circ b+b\circ a$. Thus, Corollary \ref{cor:no-deriv} implies that no Lie polynomial in a free Novikov algebra can be expressed using only symmetric Novikov products; this was previously proved by Yersaliyeva for elements of degree at most $7$ in \cite{Yersaliyeva2026}. Note that this is different from the associative case, where the double commutator $[[a,b],c]$ in an associative algebra can be expressed via the symmetrized Jordan product. 

\begin{corollary}
A non-zero linear combination of left-normed monomials 
\begin{equation}\label{eq:leftnormed}
((x_{j_1}\circ x_{j_2})\circ x_{j_3})\circ\cdots\circ x_{j_n}
\end{equation}
of degree at least $3$ in the free Novikov algebra cannot be a Lie polynomial. 
\end{corollary}

\begin{proof}
Since multilinearization preserves the properties of being a Lie polynomial and being left-normed, we may assume our Novikov polynomial to be multilinear. 
Note that the differential realization of the element \eqref{eq:leftnormed}
is $x_{j_1}x_{j_2}'\cdots x_{j_n}'$, and its symbol is therefore $t_{j_2}\cdots t_{j_n}$. A linear combination of such elements is a linear combination of elements containing every variable but one exactly once. If our element is the symbol of a Lie polynomial, it vanishes at the point 
 \[
 (0,1,\ldots,1)=(1-1,1-0,\ldots,1-0)  
 \]
which is contained in $V_{n-1}$ for $n\ge 3$, so the coefficient of $t_2\cdots t_n$ is equal to zero. Permuting the coordinates of that point, we establish that in fact all coefficients are equal to zero. 
\end{proof}

\section{A criterion for Lie polynomials in free Novikov algebras}

The main result of this section is the following criterion allowing one to recognise Lie polynomials in free Novikov algebras. 

\begin{theorem}\label{th:criterion-main}
Let $X$ be a finite set, and let $y\notin X$. Denote by $\delta_y$ the derivation of the free Novikov algebra $\Nov(X\cup\{y\})$ defined on the generators by $\delta_y(x_i)=[y,x_i]$ and $\delta_y(y)=0$. An element  $f\in\Nov(X)$ is a Lie polynomial if and only if
\begin{equation}\label{eq:NovLieCriterion}
\delta_y(f)=[y,f].    
\end{equation}
\end{theorem}

\begin{proof}
We may assume that $f$ is multilinear, for computing $\delta_y(f)$ and $[y,f]$ commutes with multilinearization. However, for multilinear elements, once we identify the Novikov operad $\Nov$ with the suboperad of $\End{W_1}$ generated by the operation $f(x)\partial\circ g(x)\partial:=f(x)g'(x)\partial$, Condition \eqref{eq:NovLieCriterion} is equivalent to the property of being $W_1$-equivariant. Indeed, an identity in the free Novikov algebra is detected by evaluating in the Novikov algebra on $\bk[x]$, and Condition \eqref{eq:NovLieCriterion} literally means equivariance with respect to the unspecified vector field corresponding to the variable $y$.
\end{proof}

This criterion is analogous to Markl's criterion allowing one to recognize Lie polynomials in pre-Lie algebras, see \cite{MR2325698}. It is demonstrated in \cite{MR4709210} that Markl's criterion may be interpreted using the theory of the so-called operadic twisting, and we shall now establish a counterpart of that result in our case. We refer the reader to \cite{MR4621635} for generalities on operadic twisting procedure; in particular, it  applies to each operad $\calP$ equipped with a morphism from the operad of Lie algebras, and we shall now describe the result in the particular case $\calP=\Nov$. Let us denote by $\alpha$ a new operation of arity~$0$ and homological degree~$-1$. The underlying non-differential operad of $\Tw(\Nov)$ is the coproduct $\Nov\vee\bk\alpha$; algebras over that operad are $\Nov$-algebras with a specified element of homological degree~$-1$ (for an expert reader who knows that the operadic twisting procedure requires to complete the coproduct with respect to the number of copies of $\alpha$, we remark that the Novikov operad is concentrated in homological degree zero, and elements with the different numbers of copies of $\alpha$ have different homological degrees, so no real completion is happening). To define the differential, one performs two steps. First, one considers the operad
 \[
\MC(\Nov)=\left(\Nov\vee\bk\alpha, d_{\MC} \right)     
 \]
encoding dg Novikov algebras with a Maurer--Cartan element; by definition, this means that the differential $d_{\MC}$ vanishes on $\Nov$ and satisfies $d_{\MC}(\alpha)=-\frac12[\alpha,\alpha]$. The element 
 \[
\ell_1^{\alpha}(a_1)=[\alpha,a_1] \in \MC(\Nov)(1)   
 \]
can be used to twist the differential of that operad, and one puts 
 \[
\Tw(\Nov)=\left(\Nov\vee\bk\alpha, d_{\Tw}=d_{\MC}+\mathrm{ad}_{\ell_1^{\alpha}}\right), 
 \]
where $\mathrm{ad}_{\ell_1^{\alpha}}(\mu)=\ell_1^{\alpha}\circ_1\mu-(-1)^{|\mu|}\sum_i\mu\circ_i\ell_1^{\alpha}$. 

\begin{corollary}\label{cor:twisting-nov}
The degree zero homology of the operad $\Tw(\Nov)$ is isomorphic to the operad $\WLie$.
\end{corollary}

\begin{proof}
Since the operad $\Tw(\Nov)$ has no elements of positive homological degrees, its degree zero homology is simply the kernel of $d_{\Tw}$ on elements of degree zero. The only elements of degree zero of $\Tw(\Nov)$ are elements of the Novikov operad. Finally, unwrapping the above definition, we see that for $f\in\Nov(n)$, we have, in notation of the previous section
 \[
d_{\Tw}(f)= [\alpha,f]-\delta_\alpha(f),    
 \] 
and hence the result follows from Theorem \ref{th:criterion-main}, since the homological degree of $\alpha$ plays no role in computing the kernel.
\end{proof}

Operadic twisting of the operad of pre-Lie algebras was used in \cite{BD24} to characterize the chain rule symmetry for singular stochastic PDEs in a vector space of an unspecified dimension. In the case of dimension one, the relevant algebraic structure collapses to Novikov algebras \cite{BD23}, and so it would be interesting to explore the consequences of Corollary \ref{cor:twisting-nov} in that context. We note that the chain rule symmetry in dimension one was successfully studied in \cite{MR5077058} by elementary methods.

\section{Open questions arising from our work}

\subsection{Invariant differential operators}

Recall that determination of the vector space 
$\Hom_{W_1}(W_1^{\otimes n}, W_1)$ 
is a particular case of the general problem of classifying invariant differential operators, that is, determining the vector space  
$\Hom_{W_1}(F_{\alpha_1}\otimes\cdots\otimes F_{\alpha_n},F_{\alpha_1+\cdots+\alpha_n+k})$ 
for $\alpha_1=\cdots=\alpha_n=-1$ and $k=n-1$. For any given values of $\alpha_1, \ldots, \alpha_n$ and $k\in\mathbb{N}$, the problem of describing that vector space seems very difficult. It is fully solved for any dimension (not just for the affine line) by Rudakov~\cite{MR0360732} for $n=1$ and by Grozman~\cite{MR575212} for $n=2$; for $n=3$, Bouarroudj~\cite{MR2547399} solved it for the affine line, and the result makes one suspect that for $n\ge 4$, the complete answer is probably not tractable. However, for given $n$ and $k$, the dimension 
 \[
\dim\Hom_{W_1}(F_{\alpha_1}\otimes\cdots\otimes F_{\alpha_n},F_{\alpha_1+\cdots+\alpha_n+k})     
 \]
is easily seen to be a semicontinuous function of $\alpha_1$, \ldots, $\alpha_n$ \cite[Sec.~III.12]{MR0463157}, that is, it is constant for Zariski generic values of these parameters, and may jump up for certain special values. The generic value of dimension in ``low degrees'' (that is, values of $k\le n$) was determined by Tabachnikov \cite[Th.~1.3]{MR742213}: for the given values of $k\le n\in\mathbb{N}$ and Zariski generic values of $\alpha_1,\ldots,\alpha_n$, we have
 \[
\dim \Hom_{W_1}(F_{\alpha_1}\otimes\cdots\otimes F_{\alpha_n},F_{\alpha_1+\cdots+\alpha_n+k})=\dim\mathcal{H}_n^k .
 \]

In view of this result, our Theorem \ref{th:LieNovInv} may be reformulated by saying that for $k=n-1$, the point $\alpha_1=\cdots=\alpha_n=-1$ is generic: the space of invariant differential operators has the same dimension as it does on a Zariski open subset. 

In fact, computer experiments and the complete classification of \cite{MR2547399} for $n=3$ suggest the following conjectures.

\begin{conjecture}\label{conj:singular values}
All values of parameters $\alpha_1,\ldots,\alpha_n$ satisfying
 \[
\sum_{i\in I} \alpha_{i} + p\ne 0,   \varnothing\ne I\subseteq\{1,\ldots,n\}, 0\le p\le |I|, 
 \]
are generic: for such values, and for $0\le k\le n-1$, we have 
 \[
\dim \Hom_{W_1}(F_{\alpha_1}\otimes\cdots\otimes F_{\alpha_n},F_{\alpha_1+\cdots+\alpha_n+k})=\dim\mathcal{H}_n^k. 
 \]
\end{conjecture}

\begin{conjecture}\label{conj:binarygeneration}
For all $0\le k\le n-1$, under the conditions on parameters described in Conjecture \ref{conj:singular values}, all elements in $\Hom_{W_1}(F_{\alpha_1}\otimes\cdots\otimes F_{\alpha_n},F_{\alpha_1+\cdots+\alpha_n+k})$ are linear combinations of elements obtained by iterated compositions of the binary operations
\begin{gather}
(f(x)dx^\alpha\otimes g(x)dx^{\beta}\mapsto f(x)g(x)dx^{\alpha+\beta})\in \Hom_{W_1}(F_{\alpha}\otimes F_{\beta},F_{\alpha+\beta}),\\
(f(x)dx^\alpha\otimes g(x)dx^{\beta}\mapsto (\alpha f(x)g'(x)-\beta f'(x)g(x))dx^{\alpha+\beta+1})\in \Hom_{W_1}(F_{\alpha}\otimes F_{\beta},F_{\alpha+\beta+1}).
\end{gather}
\end{conjecture}

Note that if we restrict ourselves to $\mathfrak{sl}_2$-invariant maps between modules of densities, then Poulain d'Andecy and Labriet \cite{MR4691919} established that generically all equivariant multilinear maps, and not just those of low degrees, can be constructed as combinations of iterations of binary ones, and moreover determined a complete set of identities between them; however, these results do not lead to any easy description of $W_1$-equivariant multilinear maps. 

\subsection{Equivariant multilinear operations on the adjoint module and Lie polynomials}

Our main theorem \ref{th:LieNovInv} can be re-stated as follows: every $W_1$-equivariant multilinear map on the adjoint module of $W_1$ is a Lie polynomial. For every Lie algebra~$\mathfrak{g}$, the Jacobi identity implies that Lie polynomials are $\mathfrak{g}$-equivariant multilinear operations on the adjoint module; however, the converse has no reason to be true. In fact, the Witt Lie algebra is the only natural example we are aware for which all equivariant multilinear operations on the adjoint module are Lie polynomials. An example where the discrepancy between equivariant multilinear operations on the adjoint module and Lie polynomials is very minor is $\mathfrak{g}=\mathfrak{sl}_2$, where we have an equality in every even arity, and an embedding of codimension one in every odd arity starting from~$3$, see \cite{MR5093214,MR3556704}. For finite-dimensional algebras, this is as good as it gets, as the following result shows.

\begin{proposition}
Let $\mathfrak{g}$ be a non-zero finite-dimensional Lie algebra over a field $\bk$ of characteristic zero. For every $n\ge 1$, there exists a $\mathfrak{g}$-equivariant operation $\mathfrak{g}^{\otimes (2n+1)}\to \mathfrak{g}$ that is not a Lie polynomial. 
\end{proposition}

\begin{proof}
We note that a non-zero finite-dimensional Lie algebra over a field $\bk$ of zero characteristic has a (possibly degenerate) non-zero invariant symmetric bilinear form $B\colon S^2(\mathfrak{g})\to\bk$. Indeed, if $\mathfrak{g}$ is perfect, it cannot coincide with its solvable radical and hence admits a non-zero semisimple quotient $\mathfrak{s}$; the pull-back of the Killing form from $\mathfrak{s}$ to $\mathfrak{g}$ is a non-zero invariant symmetric bilinear form. If $\mathfrak{g}$ is not perfect, we may take any non-zero element $\xi\in (\mathfrak{g}/[\mathfrak{g},\mathfrak{g}])^*$ and define a non-zero invariant symmetric bilinear form 
  \[
 B(x,y)=\xi(x)\xi(y).   
  \]

Once an invariant symmetric bilinear form $B$ is given, we may define the operations
 \[
x_1\otimes \cdots \otimes x_{2n+1}\mapsto B(x_1,x_2)B(x_3,x_4)\cdots B(x_{2n-1},x_{2n})x_{2n+1},      
 \]
which are equivariant and not Lie polynomials, since one can choose $x\in\mathfrak{g}$ for which $x^{\otimes(2n+1)}$ is not annihilated by this map. 
\end{proof}

It is an extremely intriguing question to determine for what other Lie algebras all equivariant multilinear maps on the adjoint module are Lie polynomials.

\printbibliography

\end{document}